\RequirePackage{fix-cm}
\documentclass[smallextended]{svjour3}       
\smartqed  
\usepackage{graphicx}
\usepackage{amssymb}
\usepackage{amsmath}
\usepackage{amsfonts}
\begin{document}

\title{Recursions on the marginals and exact computation of the normalizing constant for Gibbs processes}

\titlerunning{Marginals and normalizing constant for Gibbs processes}        

\author{C\'ecile Hardouin         \and
        Xavier Guyon
}


\institute{C. Hardouin \at
             MODALX, Universit\'e Paris Ouest Nanterre La D\'efense\\
              \email{cecile.hardouin@u-paris10.fr}
           \and
           X. Guyon \at
              SAMM, Universit\'e Paris 1
}

\date{Received: date / Accepted: date}

\maketitle

\begin{abstract}
This paper presents different recursive formulas for computing the
marginals and the normalizing constant of a Gibbs distribution $\pi.$ The
common thread is the use of the underlying Markov properties of such
processes. The procedures are illustrated with several examples, particularly
the Ising model.
\keywords{Gibbs distribution  \and  Marginal laws  \and  Normalizing constant}
\end{abstract}

\section{Introduction}
\label{intro}
Usually, the computation of the marginal distributions and the normalizing
constant $C$ of a discrete probability distribution $\pi$ involves high
dimensional summation, so that the direct evaluation of these sums becomes
quickly infeasible in practice. For example, for an Ising model on a
$10\times10$ grid, it involves a summation over $2^{100}$ terms. This problem
has a deep impact in many applications, as for instance maximum likelihood
 estimation, and some significant efforts have been undertaken to bypass the
problem. In spatial statistics, we replace the likelihood by the conditional
pseudo likelihood (Besag 1974; Gaetan and Guyon 2010). Another way is to
estimate $C$ using efficient Monte Carlo methods, see for example Moeller et
al. (2006) or Huber (2012). However, it is sometimes possible to compute $C$
exactly using an efficient algorithm; a first description of such an algorithm
is given in Liu (2001); Pettitt et al. (2003) obtain the exact normalizing
constant for a general categorical $N-$ valued distribution on a $m\times n$
cylinder with a method involving the computation of the eigenvalues of
a $N^{m}\times N^{m}$ matrix, which makes it possible for $N^{m}\lesssim1024$,
that is for example $N=2$ and $m=10.$ Reeves and Pettitt (2004) give
recursions for a factorizable distribution $\pi$. Then Friel and Rue (2007)
extend Reeves and Pettitt's results to allow exact sampling, maximization of
distributions and computation of marginal distributions. We propose here other
procedures to compute the marginals and the normalizing constant for Gibbs
processes. We express the Gibbs distribution in a matrix form; our
first result then yields an expression for the normalizing constant; this result
coincides with that of Reeves and Pettitt (2004), however we formulate it
differently, as a simple matrix product, which facilitates fast computation; indeed, new techniques emerging from general big data analysis can be applied here; on the one hand, matrix-vector and matrix-matrix multiplication is efficient in many architectures, for instance GPUs, multi-core processors, or parallel computers.  On the other hand, eigenvalue decomposition or low-rank matrix factorizations for large matrices can be performed using efficient tools, such as probabilistic algorithms (Halko 2011, Witten 2013).
We also give the marginals on a
general subset. We then give a new algorithm to compute the marginal distributions, based
on future conditional probabilities. Finally, we propose a new method to
compute dichotomous marginals, based on a Markov chain - Markov random field framework.

The paper is organized as follows. We first consider a temporal random
variable $Z(T)=(Z_{1},Z_{2},\cdots,Z_{T})$ on a finite state space, with a
Gibbs distribution $\pi$. Section 2 summarizes basic properties about this
process and we then present results which allow the computation of the
marginals and the normalization constant of $\pi$ under a matrix formulation. We
also compare the numerical efficiency of various algorithms in terms of their
computing times. Section 3 presents two new algorithms computing
marginal distributions of $\pi$; the first one is based on the conditionals on the future and
involves the normalizing constant, while the second one recursively computes
the distribution $\pi$ as well as its dichotomous sequence of marginals, avoiding the
normalizing constant problem; this last procedure directly follows from the
Markov properties of the process. We extend the results in section 4 for
general Gibbs fields, illustrating the result with an Ising model on a
lattice; some numerical experiments are given. The paper ends with some generalisations and opens a discussion about large lattices; the matrix formulation combined with randomized algorithms for low-rank matrix approximations, and new architectures, lead to efficient computation methods.

\section{Recursions for a temporal Gibbs distribution}

\label{sec:1}

\subsection{Factoring joint distribution, Markov field and Markov chain
properties}

\label{factorisante 1} Let $T>0$ be a positive integer, $\mathcal{T}%
=\{1,2,\cdots,T\},~E$ be a finite state space with $N$ elements, $Z(T)=(Z_{1}%
,Z_{2},\cdots,Z_{T})$ be a temporal random variable with a factorizable joint
distribution $\pi$ on $E^{T}$ (see Reeves and Pettitt 2004), that is,
\begin{equation}
\pi(z_{1},\cdots,z_{T})=C^{-1}\exp\sum_{s=1,T-1}h_{s}(z_{s},z_{s+1}%
)\label{gibbs}%
\end{equation}

Let us denote by $\pi_{A}$ the marginal of $\pi$ on the subset $A\subseteq
\mathcal{T}$,  $\pi_{s}^{t}=\pi_{\{s,s+1,\cdots,t\}}$ with $s<t$, and
$u_{s}^{t}=(u_{s},...,u_{t})\in E^{t-s+1}.$

We can interpret $\pi$ as a Gibbs distribution with energy $U_{T}(z_{1}%
,z_{2},..z_{T})=\sum_{s=1,T-1}h_{s}(z_{s},z_{s+1})$ associated with the saturated
pair potentials $(h_{s})_{s=1,T-1}$. In particular, $\pi$ is a bilateral
Markov random field w.r.t. the 2-nearest neighbors system (see Kindermann et
Snell 1980; Lauritzen 1996; Guyon 1995):%
\[
\pi(z_{t}\mid z_{s},s\neq t)=\frac{\exp\left(  h_{t-1}(z_{t-1},z_{t}%
)+h_{t}(z_{t},z_{t+1})\right)  }{\sum_{u\in E}\exp\left(  h_{t-1}%
(z_{t-1},u)+h_{t}(u,z_{t+1})\right)  }=\pi(z_{t}\mid z_{t-1},z_{t+1}).
\]

\noindent This non-causal conditional distribution $\pi(z_{t}\mid
z_{t-1},z_{t+1})$ can be easily computed. On the other hand, the normalizing
constant $C$ as well as the marginals generally cannot be computed since they entail
high dimensional summation.

\medskip

\noindent Note that $Z(T)$ is also a Markov chain; indeed, $\pi(z_{t}\mid z_{s},~s\leq
t-1)=\pi_{1}^{t}(z_{1}^{t})  / \pi_{1}^{t-1}(z_{1}^{t-1})$; let us  write $\pi_{1}^{t}(z_{1}^{t})=\sum\limits_{u_{t+1}%
^{T}\in E^{T-t}}\pi(z_{1},\cdots,z_{t},u_{t+1},\cdots,u_{T})$   for $1\leq t\leq T$.  Finally, we
have
\[%
\begin{tabular}
[c]{ll}%
$\pi(z_{t}\mid z_{s},~s\leq t-1)$ & $=\frac{\exp h_{t-1}(z_{t-1},z_{t}%
)\times\sum\limits_{u_{t+1}^{T}}\exp\left(  h_{t}(z_{t},u_{t+1})+\sum
_{s=t+1}^{T-1}h_{s}(u_{s},u_{s+1})\right)  }{\sum\limits_{u_{t}^{T}}%
\exp\left(  h_{t-1}(z_{t-1},u_{t})+\sum_{s=t}^{T-1}h_{s}(u_{s},u_{s+1}%
)\right)  }$\\
& $=\pi(z_{t}\mid z_{t-1}).$%
\end{tabular}
\]
Generally the analytic form of this unilateral transition is not explicit due to
the high dimensional summation in $u_{t+1}^{T}$. However, an explicit computation might be sometimes
possible, in which case we will present in \ref{dich} a method for computing a sequence of dichotomous marginals, as well as the full distribution, avoiding
the normalizing constant computation.

\subsection{Forward recursions over marginal distributions}

We first introduce the matrix formulation which will be used throughout the
paper. Let us define $H_{s}$, $s=1,T-1$, to be  $N\times N$ matrices with components
$H_{s}(u,v)=\exp h_{s}(u,v)$. For such two $N\times N$ matrices $H$ and $G$,
for a $N$-row vector $F=(F(v),v\in E),$ and for a $N$-column vector
$B=(B(v),v\in E)$, we introduce, in a classical way, the matrix products
$HG(u,v)=\displaystyle \sum_{w\in E}H(u,w)G(w,v)$, $\displaystyle FH(v)=\sum_{u\in
E}F(u)H(u,v)$ and $HB(u)=\displaystyle \sum_{v\in E}H(u,v)B(v)$.

Coming back to (\ref{gibbs}), we write
\[
\pi(z_{1},...,z_{T})=C^{-1}\prod_{s=1,T-1}H_{s}(z_{s},z_{s+1}).
\]

\medskip

Now, let  $B_{T}=\mathbf{1}$ be the $N$-column vector with constant
coordinates $1$. A direct marginalization on $z_{T}$ gives:%
\[
\pi_{1}^{T-1}(z_{1},...,z_{T-1})=C^{-1}\left \{\prod_{s=1,T-2}H_{s}(z_{s}%
,z_{s+1})\right\}H_{T-1}B_{T}(z_{T-1}).
\]
In the same way, defining $B_{t-1}=H_{t-1}B_{t}$, we have for $t=T,2$ :%
\begin{align}
\pi_{1}^{t}(z_{1},...,z_{t}) &  =C^{-1}\prod_{s=1,t-1}H_{s}(z_{s}%
,z_{s+1})(H_{t}B_{t+1})(z_{t})\label{recursion arriere}\\
&  =C^{-1}\prod_{s=1,t-1}H_{s}(z_{s},z_{s+1})(H_{t}\cdots H_{T-1}B_{T}%
)(z_{t})~.\nonumber
\end{align}
These equations, up to the constant $C$, give the marginals $\pi_{1}^{t}$.
Moreover, looking at $\pi_{\{1\}}$, we obtain the expression of the normalizing constant.

\begin{proposition}
1- The normalizing constant of the general model (\ref{gibbs}) is given by
\begin{equation}
C=B_{T}^{\mathrm{T}}\left\{\prod
_{s=1,T-1}H_{s}\right\}B_{T}.\label{calcul de la constante}%
\end{equation}

2- For time invariant potentials , the formula reduces to $C=\mathbf{1}%
^{\mathrm{T}}(H)^{T-1}\mathbf{1.}$
\end{proposition}

\textbf{Proof:}

Indeed, from (\ref{recursion arriere}), equation $C=C\times\sum_{z_{1}\in E}%
\pi_{\{1\}}(z_{1})$  gives the first  result. Then, for the time
invariant potentials $h_{s}=h$,  $H_{s}=H$  implies
$C=\mathbf{1}^{\mathrm{T}}(H)^{T-1}\mathbf{1,}$ or $C=\mathbf{1}^{\mathrm{T}%
}(H)^{T-2}H_{T-1}\mathbf{1}$ if $H_{T-1}\neq H$.

\bigskip

These results agree with those of Liu (2001,
paragraph 2.4), and Reeves and Pettitt (2004). However, writing
(\ref{calcul de la constante}) in terms of matrix products opens the possibility of other
implementation methods; practically, as we will see in the next examples, we can sometimes improve the computing times if we diagonalize the
matrix $H$.

Let us now give our results on the computation of marginal laws on general subsets. Still using
forward recursions, and defining $F_{t}=F_{t-1}H_{t-1}$ for $t\geq2$, we
obtain :
\[
\pi_{t}^{T}(z_{1},...,z_{T})=C^{-1}\left\{(F_{t-1}H_{t-1})(z_{t})\prod
_{s=t,T-1}H_{s}(z_{s},z_{s+1})\right\}
\]%
\[
=C^{-1}\left\{(F_{1}H_{1}H_{2}\cdots H_{t-1})(z_{t})\prod_{s=t,T-1}H_{s}%
(z_{s},z_{s+1})\right\}.
\]

More generally, we obtain the following result:

\begin{proposition}
1 - The marginal distribution of $\pi$ on $S=\{s_{1},s_{2},\cdots
,s_{q}\}\subseteq\mathcal{T}$ with $1=s_{1}<s_{2}<\cdots<s_{q-1}<s_{q}=T$ is
given by
\begin{equation}
\pi_{S}(z_{S})=C^{-1}\prod_{i=1,q-1}\left (\prod_{s=s_{i}}^{s_{i+1}-1}%
H_{s}\right)(z_{s_{i}},z_{s_{i+1}}),\label{marginale generale}%
\end{equation}
2 - The marginal on $S_{1,T}=S\setminus\{1,T\}$ is obtained by replacing the
first $H$-product $\displaystyle\left(\prod_{1}^{s_{2}-1}H_{s}%
\right)(z_{1},z_{s_{2}})$ with $\mathbf{1}^{\mathrm{T}}\displaystyle\left(\prod
_{1}^{s_{2}-1}H_{s}\right)(z_{s_{2}})$ and the last product $\displaystyle\left(\prod
_{s_{q-1}}^{s_{q}-1}H_{s}\right)(z_{s_{q-1}},z_{T})$ with $\left(\left(\displaystyle\prod
_{s_{q-1}}^{s_{q}-1}H_{s}\right)\mathbf{1}\right)(z_{s_{q-1}})$.
\end{proposition}

\subsection{Examples and numerical performances}

 Coming back to formula (\ref{calcul de la constante}), if the size $N$ of $E$ allows the
diagonalization of the matrix $H$, we can achieve the computation of $C$
independently of the temporal dimension $T$. Let us look at two examples for
which we compare the computing times for different algorithms. For this study,
we have used the Matlab software (programs are available on request).

\bigskip

\textit{Example 1 : binary temporal model}

 The following toy example illustrates the relative efficiency of different methods.  Let us consider $E=\{0,1\}$ and the auto-logistic Gibbs field $\pi$
associated with time independent singleton and pair potentials $\theta
_{t}(z_{t})=\alpha z_{t}$, $t=1,T$ and $\Psi_{t}(z_{t},z_{t+1})=\beta
z_{t}z_{t+1}$ for $t\leq T-1$. We have $h_{t}(z_{t},z_{t+1})=\theta_{t}%
(z_{t})+\Psi_{t}(z_{t},z_{t+1})$ for $t=1,T-2$ and $h_{T-1}(z_{T-1}%
,z_{T})=\theta_{T-1}(z_{T-1})+\Psi_{T-1}(z_{T-1},z_{T})+\theta_{T}(z_{T}).$

We present in Table \ref{tab:1} the times for the computation of $C$ for increasing
values of $T$ and the following methods:

Method 1: $C=$ $^{t}\mathbf{1}(H)^{T-2}H_{T-1}\mathbf{1}$;

Method 2: \textit{\ }$H$ is diagonalized, with the Matlab function ``eig''.

Method 2bis: $H$ is diagonalized, with the diagonal and
eigenvectors matrices computed analytically by hand.

Method 3: Complete computation of $C=\left[  (\lambda
_{1}^{T-1}(\lambda_{2}-1)+\lambda_{2}^{T-1}(1-\lambda_{1}))(1+e^{\alpha
}) \right. $

\noindent $\left. +e^{\alpha}(1+e^{\alpha+\beta})(\lambda_{2}^{T-1}-\lambda_{1}%
^{T-1})\right]  /\sqrt{\Delta} $,  using the eigenvalues $\lambda_{1}=(1+e^{\alpha+\beta
}-\sqrt{\Delta})/2$ and $\lambda_{2}=(1+e^{\alpha+\beta}+\sqrt{\Delta})%
/2$ , with $\Delta=(1-e^{\alpha+\beta})^{2}+4e^{\alpha}$.

Method 4: $C$ is obtained by summation, using a bitmap dodge which computes
simultaneously the $2^{T}$ elements of $E^{T}$. Practically we use Matlab functions dec2bin and str2num which respectively convert decimals to binary numbers in string and conversely  strings to numbers.

Method 4bis: $C$ is obtained by summation, with a simple loop (each element is
computed individually and added to the previous computation).

We stopped computing $C$ by summation (methods 4 and 4 bis) for $T>25.$

\begin{table}
\caption{Computing times (in seconds) of $C$\ for a binary temporal Gibbs
distribution, $\alpha=1$, $\beta=-0.8$.}
\label{tab:1}       
\begin{tabular}
[c]{|l|c|c|c|c|c|c|c|}\hline
& \textit{Meth. 1} & \textit{Meth. 2} & \textit{Meth.~2 bis} & \textit{Meth.
3} & \textit{Meth. 4} & \textit{Meth. 4 bis} & Value of $C$\\\hline
$T=10$ & \multicolumn{1}{|c|}{$8.2\times10^{-5}$} &
\multicolumn{1}{|c|}{$1.1\times10^{-4}$} & \multicolumn{1}{|l|}{$7.1\times
10^{-5}$} & \multicolumn{1}{|c|}{$1.8\times10^{-5}$} &
\multicolumn{1}{|l|}{$0.1134$} & $0.2150$ & $3.3441\times10^{4}$\\\hline
$T=20$ & \multicolumn{1}{|c|}{$9.2\times10^{-5}$} &
\multicolumn{1}{|c|}{$1.2\times10^{-4}$} & \multicolumn{1}{|l|}{$7.9\times
10^{-5}$} & \multicolumn{1}{|c|}{$2.0\times10^{-5}$} &
\multicolumn{1}{|l|}{$204.62$} & $375.08$ & $8.6756\times10^{8}$\\\hline
$T=25$ & \multicolumn{1}{|c|}{$1.3\times10^{-4}$} &
\multicolumn{1}{|c|}{$1.9\times10^{-4}$} & \multicolumn{1}{|l|}{$5.2\times
10^{-5}$} & \multicolumn{1}{|c|}{$1.9\times10^{-5}$} &
\multicolumn{1}{|l|}{$8118.2$} & $\sim3.9$ hours &
\multicolumn{1}{|l|}{$1.3974\times10^{11}$}\\\hline
$T=500$ & \multicolumn{1}{|c|}{$1.1\times10^{-4}$} &
\multicolumn{1}{|l|}{$1.3\times10^{-4}$} & \multicolumn{1}{|l|}{$7.9\times
10^{-5}$} & \multicolumn{1}{|l|}{$2.4\times10^{-5}$} & \multicolumn{1}{|l|}{}
&  & \multicolumn{1}{|l|}{$6.4759\times10^{220}$}\\\hline
\end{tabular}
\end{table}

We observe that the computing times of $C=$ $^{t}\mathbf{1}%
(H)^{T-2}H_{T-1}\mathbf{1}$ are negligible for $T<500$ for methods 1 to 3, while methods 4 and
4bis become quickly unusable. There is a slight advantage to diagonalization
by hand, since Matlab diagonalization involves a (small) computing time; on the
other hand, diagonalization is not worthwhile since the matrix's size is only
$2\times2$. We present in Table \ref{tab:2} the computation of $C$ for $T>600$, performed with    method 1. We observe that the computing time does not increase with $T$, it remains stable (and small). 

\begin{table}
\caption{Computing times (in seconds) of $C$ for a binary temporal Gibbs
distribution, $\alpha=1$, $\beta=-0.8$.}
\label{tab:2}       
\begin{tabular}
[c]{|l|c|c|c|}\hline
& $T=1000$ & $T=10000$ & $T=10^{6}$\\\hline
\textit{Method 1} & \multicolumn{1}{|c|}{$8.1\times10^{-5}$} & $6.2\times
10^{-5}$ & $7.7\times10^{-5}$\\\hline
Value of  $C$ & \multicolumn{1}{|c|}{$9.9392\times2^{500}\times10^{290}$} &
$1.9918\times2^{900\times16}\times10^{79}$ & $1.4010\times2^{880\times
1666}\times10^{68}$\\\hline
\end{tabular}
\end{table}

\medskip
\textit{Example 2 : bivariate binary temporal model}

We now consider  $E=\{0,1\}^{2}$ and $Z(T)$ an anisotropic Ising
model with invariant saturated potentials $h_{s}=h$,
$$
h((x_{1},y_{1}),(x_{2},y_{2}))=\alpha x_{1}+\beta y_{1}+\gamma x_{1}%
y_{1}+\alpha x_{2}+\beta y_{2}+\gamma x_{2}y_{2}+\delta(x_{1}x_{2}+y_{1}%
y_{2}),$$
with the convention that a pair potential equals $0$ if a state is taken out
of the time domain $\{1,2,\cdots,T\}$. We computed the constant $C$ in two
ways, first directly calculating the power $H^{T-2}$, then making use of the
diagonalization of $H,$ i.e. calculating $C=$ $^{t}\mathbf{1}PD^{T-2}%
P^{-1}H_{T-1}\mathbf{1}$. The parameter values are $\alpha=1,~\beta
=-0.8,~\gamma=-0.5,~$ and $\delta=0.04$.  Table 3 displays some results; we
see that since we now deal with a matrix of size $4\times4,$ diagonalization is
more attractive than in the previous example.

\begin{table}
\caption{Computing times (in seconds) of $C$\ for a bivariate binary temporal
Gibbs distribution.}
\label{tab:3}       
\begin{tabular}
[c]{|l|c|c|c|}\hline
& Method 1 & Method 2 & Value $C$\\\hline
$T=500$ & $0.0074$ & $3.8\times10^{-4}$ & $9.1687\times2^{400}\times10^{229}%
$\\\hline
$T=1000$ & $0.0059$ & $3.5\times10^{-4}$ & $2.0726\times2^{1400}\times
10^{279}$\\\hline
$T=10000$ & $0.0079$ & $4.6\times10^{-4}$ & $4.8390\times2^{930\times24}%
\times10^{288}$\\\hline
$T=10^{6}$ & $0.0130$ & $0.0043$ & $8.2420\times2^{930\times2500+2670}%
\times10^{69}$\\\hline
\end{tabular}
\end{table}

\subsection{Results for $r$-range potentials}

Reeves and Pettitt (2004) consider more general $r$-factorizable distributions
$\pi(z(T))=C^{-1}\prod_{s=1}^{T-r}H_{s}(z_{s},z_{s+1},\cdots,z_{s+r})$. There,
the function $H_{s}$ is defined on $E^{\ast}=E^{r+1}$. For a real function $H$
defined on $E^{\ast}$, we define $H^{\ast}$ on $E^{\ast}\times E^{\ast}$
as:%
\begin{equation}
H^{\ast}(u,v)=H(u_{2},...,u_{r+1},v_{r+2})\prod_{i=1}^{r}\mathbf{1}%
(u_{i+1}=v_{i}).\label{H etoile}%
\end{equation}
Then, with respect to the $(\ast)$-objects, and with the notations of
\S \ref{factorisante 1}, we obtain the same results; for instance the
normalizing constant equals $C=\mathbf{1}^{\mathrm{T}}(\prod_{s=1}^{T-r}%
H_{s}^{\ast})\mathbf{1}$. Recursive algorithms for the marginals of $\pi$
follow in the same way as in (\ref{recursion arriere}) and
(\ref{marginale generale}).

\section{Other recursions}

\subsection{A new recursive algorithm for marginals based on future
conditionals}

For the sake of simplicity, let us specify the Gibbs specification
(\ref{gibbs}) of $\pi$ in terms of singleton and pair potentials, and write
:%
\[
h_{s}(z_{s},z_{s+1})=\theta_{s}(z_{s})+\Psi_{s}(z_{s},z_{s+1})\ \mathrm{\ for}%
\ s=1,T\ ,
\]
with the convention $\Psi_{T}\equiv0$. We then get the following expression for the energy for
$t=1,T$:%
\[
U_{t}(z_{1},..z_{t})=\sum_{s=1,t}\theta_{s}(z_{s})+\sum_{s=1,t-1}\Psi
_{s}(z_{s},z_{s+1})
\]

We compute recursively the marginal $\pi_{1}^{t}$ by conditioning on the
future. For $t<T$, it is clear that $\pi(z_{1},..z_{t}\mid z_{t+1}%
,..z_{T})=\pi(z_{1},..z_{t}\mid z_{t+1})$ with
\[
\pi(z_{1},..z_{t}\mid z_{t+1})=C_{t}^{-1}(z_{t+1})\exp U_{t}^{\ast}%
(z_{1},..z_{t};z_{t+1}),
\]
where
\begin{equation}
U_{t}^{\ast}(z_{1},..z_{t};z_{t+1})=U_{t}(z_{1},..z_{t})+\Psi_{t}%
(z_{t},z_{t+1})\label{future cond energy}%
\end{equation}
is the \textit{future-conditional} energy, and $C_{t+1}(z_{t+1})=\sum
_{u_{1}^{t}\in E^{t}}\exp\left\{  U_{t}^{\ast}(z_{1}^{t};z_{t+1})\right\}  $.

\bigskip

For $1\leq t\leq T,$ we define $\gamma_{t}(z_{1},z_{2},\cdots,z_{t};u)=\exp
U_{t}^{\ast}(z_{1}^{t};u)$ and the $N-$row vector
\[
\Gamma_{t}(z_{1},z_{2},\cdots,z_{t})(u)=\gamma_{t}(z_{1},z_{2},\cdots
,z_{t};u)~,~~\ ~u\in E~.
\]

For $t=T$, $\Gamma_{T}(z(T))$ is the constant vector with components $\gamma
_{T}(z(T))=\exp U_{T}(z(T))$. Notice that $\Gamma_{t}(z_{1},..z_{t})$ is
analytically explicit.

For $1\leq t\leq T$ and $u,v\in E$, $H_{t}(u,v)=\exp\{\theta_{t}(u)+\Psi
_{t}(u,v)\}$; let us define the sequence $(D_{t},t=T,2)$ of $N$-column vectors
by $D_{T}=$ $^{t}(1,0,\cdots,0),$ and $D_{t-1}=H_{t}D_{t}$ for $t\leq T$. The
 following result then gives a new recursion for the marginals.

\begin{proposition}
\label{new recursion}Recursion for marginal distributions.

1 - For $2\leq t\leq T$, we have:
\begin{equation}
\sum_{z_{t}\in E}\Gamma_{t}(z_{1},z_{2},\cdots,z_{t})=\Gamma_{t-1}(z_{1}%
,z_{2},\cdots,z_{t-1})H_{t}.\label{recurrence fondamentale}%
\end{equation}

2 - For $1\leq t\leq T$ ,%
\begin{equation}
\pi_{1}^{t}(z_{1},z_{2},\cdots,z_{t})=C^{-1}\times\Gamma_{t}(z_{1}%
,z_{2},\cdots,z_{t})D_{t}.\label{lois marginales}%
\end{equation}

3 - These equations give a new expression for the normalizing constant:%
\begin{equation}
C=D_{T}^{\mathrm{T}}\left\{\prod_{s=2,T}H_{s}\right\}D_{T}.
\end{equation}

\end{proposition}

\textbf{Proof:}

1 - For $2\leq t\leq T$, $U_{t}(z_{1},..z_{t-1},z_{t})=U_{t-1}(z_{1}%
,..z_{t-1})+\theta_{t}(z_{t})+\Psi_{t-1}(z_{t-1},z_{t})$; therefore,
\begin{align*}
U_{t}^{\ast}((z_{1},..z_{t-1},z_{t});z_{t+1}) & =U_{t-1}(z_{1},..z_{t-1}%
)+\theta_{t}(z_{t})+\Psi_{t-1}(z_{t-1},z_{t})+\Psi_{t}(z_{t},z_{t+1})\\
& =U_{t-1}^{\ast}(z_{1},..z_{t-1};z_{t})+\left[  \theta_{t}(z_{t})+\Psi
_{t}(z_{t},z_{t+1})\right]  .
\end{align*}

\noindent This implies $\gamma_{t}((z_{1},..z_{t-1},z_{t});u)=\gamma
_{t-1}(z_{1},..z_{t-1};z_{t})\times H_{t}(z_{t},u)$ and summation over $z_{t}
$ gives (\ref{recurrence fondamentale}).

\medskip2 - We prove (\ref{lois marginales}) by step-down induction. For
$t=T$, the equality holds since $\pi_{1}^{T}(z_{1},..z_{T})=\pi
(z_{1},..z_{T})=C^{-1}\exp U_{T}(z_{1},..z_{T})=C^{-1}\times\Gamma_{T}%
(z_{1},..z_{T})D_{T}.$

Let us assume that (\ref{lois marginales}) is satisfied for some $t,$ $2\leq
t\leq T$. We use (\ref{recurrence fondamentale}) which gives:%
\begin{align*}
\pi_{1}^{t-1}(z_{1}^{t-1})  & =\sum_{z_{t}}\pi_{1}^{t}(z_{1}^{t-1}%
,z_{t})=C^{-1}\left\{\sum_{z_{t}}\Gamma_{t}(z_{1},..z_{t-1},z_{t})\right\}D_{t}\\
& =C^{-1}\Gamma_{t-1}(z_{1},..z_{t-1})H_{t}D_{t}=C^{-1}\Gamma_{t-1}%
(z_{1},..z_{t-1})D_{t-1}.
\end{align*}
This completes the proof.

\bigskip

Proposition (\ref{new recursion}) can be extended in a natural way to $r$-lag
Gibbs processes. For example, let us consider the $2$-lag factorizable
distribution $\pi$, characterized by the energy: $U_{T}(z_{1},..z_{T}%
)=\displaystyle\sum_{s=1,T}\theta_{s}(z_{s})+\sum_{s=1,T-1}\Psi_{1,s}%
(z_{s},z_{s+1})+\sum_{s=1,T-2}\Psi_{2,s}(z_{s},z_{s+2})$ with the convention
$\Psi_{1,T}\equiv\Psi_{2,T-1}\equiv\Psi_{2,T}=0$.
 It is easy to see that $\pi$ is a Markov field w.r.t. the 4-nearest neighbors
system and $\pi(z_{1},..z_{t}\mid z_{t+1},..z_{T})=\pi(z_{1},..z_{t}\mid
z_{t+1},z_{t+2})=C_{t}^{\ast}(z_{t+1},z_{t+2})\exp U_{t}^{\ast}(z_{1}%
^{t};z_{t+1},z_{t+2}),$ where the future conditional energy equals
$U_{t}^{\ast}(z_{1}^{t};z_{t+1},z_{t+2})=U_{t}(z_{1}^{t})+\Psi_{1,t}%
(z_{t},z_{t+1})+\Psi_{2,t-1}(z_{t-1},z_{t+1})+\Psi_{2,t}(z_{t},z_{t+2}).$

Again $U_{t}^{\ast}(z_{1}^{t};z_{t+1},z_{t+2})=U_{t-1}^{\ast}(z_{1}%
^{t-1};z_{t},z_{t+1})+\theta_{t}(z_{t})+\Psi_{1,t}(z_{t},z_{t+1})+\Psi
_{2,t}(z_{t},z_{t+2}).$

\noindent Following the previous scheme, we define for $t\leq T,$ the $N^{2}%
$-row vector $\Gamma_{t}(z_{1},..z_{t})$ as
\[
\Gamma_{t}(z_{1},..z_{t})(u,v)=\exp U_{t}^{\ast}(z_{1},..z_{t};u,v),~~~u,v\in
E~.
\]
Then, as in the proof of Proposition \ref{new recursion}, but w.r.t the
matrices $H_{s}^{\ast}$ (\ref{H etoile}) defined on $E^{\ast}\times E^{\ast}$, where
$E^{\ast}=E^{2}$, we obtain the recurrence
(\ref{recurrence fondamentale}) on the contributions $\Gamma_{t}(z_{1}%
,..z_{t})$ and the result (\ref{lois marginales}) on the marginals.

\subsection{Dichotomous marginals sequence\label{dich}}

Let us consider the following dichotomous sequence of embedded subsets of $\mathcal{T}%
=\{1,2,...,T\}.$ We set $T=2^{r}+1$ and introduce the sequence $S=S_{r+1}%
=\{1,2,...,T\},~S_{r}=\{1,3,5,...,T\},~S_{r-1}=\{1,5,9,...T\},\dots $
,$S_{1}=\{1,T\}.$

\bigskip

Let us denote by $\pi_{r+1}=\pi$ the distribution of $Z(T)=(Z_{1},Z_{2}%
,...,Z_{T})$ on $S=S_{r+1}.$ We recall that $\pi$ is a bilateral Markov random
field w.r.t. the 2-nearest neighbors system: $\pi(z_{t}\mid z_{s},~s\neq
t)=\pi(z_{t}\mid z_{t-1},z_{t+1}).$ On the other hand, $Z(T)$ is also a Markov
chain; let $P$ be its transition matrix. The following expression gives the
links between the Markov field parameters and components of the transition matrix:%

\begin{equation}
\pi(z_{t}\mid z_{t-1},z_{t+1})=\frac{\pi(z_{t+1}\mid z_{t})\pi(z_{t}\mid
z_{t-1})}{\sum\limits_{u}\pi(z_{t+1}\mid u)\pi(u\mid z_{t-1})}%
~.\label{chaine-champ}%
\end{equation}

We assume here that we can find this explicit link (as for instance in the
example below). Let us keep every other term of $Z$; we get a process
$(Z_{1},Z_{3},...,Z_{T})$ which is a Markov chain with transition matrix
$P^{2};$ therefore, $(Z_{1},Z_{3},...,Z_{T})$ is also a bilateral Markov
random field with distribution $\pi_{r}$\ linked to $P^{2}$ by an
equation analogous to (\ref{chaine-champ}), with conditional distributions $\pi
_{r}(z_{t}\mid z_{t-2},z_{t+2}).$ We iterate the process: $(Z_{1}%
,Z_{5},..Z_{T})$ is a Markov chain with transition matrix $P^{4}$ and a
bilateral Markov random field with distribution $\pi_{r-1}$. The process is repeated until
we get $(Z_{1},Z_{T}).$

We write $\pi(z_{1},z_{2},...,z_{T})=\pi(z_{1},z_{2},...,z_{T}\mid z_{1}%
,z_{3},z_{5},..,z_{T})\times\pi(z_{1},z_{3},z_{5},..,z_{T}).$ Next, on the one
hand,
$$\pi(z_{1},z_{2},...,z_{T}\mid z_{1},z_{3},z_{5},..,z_{T})=\pi(z_{2}\mid
z_{1},z_{3})\times\pi(z_{4}\mid z_{3},z_{5})...\times\pi(z_{T-1}\mid z_{T-2},z_{T}
)$$

$=\prod\limits_{s\in S\setminus S_{r}}\pi(z_{s}\mid z_{s-1},z_{s+1}).$

On the
other hand, $\pi(z_{1},z_{3},z_{5},..,z_{T})=\pi_{r}(z_{1},z_{3},z_{5}%
,z_{7},...,z_{T}\mid z_{1},z_{5},z_{9},..,z_{T})\times\pi_{r}(z_{1}%
,z_{5},z_{9},..,z_{T})$ which equals as previously $\prod\limits_{s\in
S_{r}\setminus S_{r-1}}\pi_{r}(z_{s}\mid z_{s-1},z_{s+1})\times\pi_{r}%
(z_{1},z_{5},z_{9},..,z_{T})$.

Repeating the process, we obtain
\[
\pi(z_{1},z_{2},...,z_{T})=\prod\limits_{s\in S\setminus S_{r}}\pi(z_{s}\mid
z_{s-1},z_{s+1})\times\prod\limits_{s\in S_{r}\setminus S_{r-1}}\pi_{r}%
(z_{s}\mid z_{s-1},z_{s+1})\times...
\]
\begin{equation}
\times\prod\limits_{s\in S_{2}\setminus S_{1}}\pi_{2}(z_{s}\mid z_{s-1}%
,z_{s+1})\times\pi_{1}(z_{1},z_{T})~.\label{pi-final}%
\end{equation}

This gives a new algorithm for computing $\pi(z_{1},z_{2},...,z_{T})$ or
any \textquotedblleft dichotomous\textquotedblright\ marginal on a subset
$S_{j};~$\ we also derive the normalizing constant with $\pi(z_{1}%
,z_{2},...,z_{T})=C^{-1}\exp\sum\limits_{s=1,T-1}h_{s}(z_{s},z_{s+1}).$

\medskip

\textit{Example 3 :  the classical Ising model}

We consider the state space $E=\{-1,+1\},$ and the distribution $\pi$\
defined by $h_{s}(z_{s},z_{s+1})=\alpha z_{s}+\beta z_{s}z_{s+1}.\ $

We also write $Z$'s transition matrix as $P=\left(
\begin{array}
[c]{cc}%
p & 1-p\\
1-q & q
\end{array}
\right)  .$ Equation (\ref{chaine-champ}) amounts to%
\[
\frac{\exp\left(  \alpha z_{t}+\beta(z_{t-1}z_{t}+z_{t}z_{t+1})\right)  }%
{\exp\left(  \alpha+\beta(z_{t-1}+z_{t+1})\right)  +\exp-\left(  \alpha
+\beta(z_{t-1}+z_{t+1})\right)  }=\frac{\pi(z_{t+1}\mid z_{t})\pi(z_{t}\mid
z_{t-1})}{\sum\limits_{u}\pi(z_{t+1}\mid u)\pi(u\mid z_{t-1})}%
\]
and it is easy to see that links between the Markov field parameters and the
transition matrix are given by the following formula (see also Guyon 1995)
\begin{equation}
\alpha=\frac{1}{2}\ln\frac{p}{q}\text{ and }\beta=\frac{1}{4}\ln\frac
{pq}{(1-p)(1-q)}\label{alphabetapq}%
\end{equation}

Let us consider the following step. Keeping every other term of $Z$ yields a Markov
chain with transition $P^{2}$ and a Markov field with parameters $\alpha_{2} $
and $\beta_{2}$ linked to the elements of $P^{2}$ by equations analogous to
(\ref{alphabetapq}). Some further calculations then give $\alpha_{r}$ and $\beta_{r}$
which depend on $\alpha$ and $\beta.$

Finally, we set for $j=r$ down to $1$ the following equation
\begin{align*}
\alpha_{j}  &  =\frac{1}{2}\ln\left[  \frac{1+e^{2\alpha_{j+1}+4\beta_{j+1}}%
}{1+e^{-2\alpha_{j+1}+4\beta_{j+1}}}\right] \\
\beta_{j}  &  =\frac{1}{4}\ln\left[  1+\frac{e^{4\beta_{j+1}-2\alpha_{j+1}%
}(1-e^{-4\beta_{j+1}})^{2}}{(1+e^{-2\alpha_{j+1}})^{2}}\right]
\end{align*}

with $\alpha_{r+1}=\alpha$ and $\beta_{r+1}=\beta.$

We recursively compute all parameters $\alpha_{j}$ and $\beta_{j}$ which
allows the computation of $\pi$ using (\ref{pi-final}).

\bigskip

This method allows for a  true maximum likelihood estimation, and does not
require a normalizing constant computation.
\section{The case of spatial Gibbs fields}

\subsection{A temporal multidimensional Gibbs process}

Let us consider $Z=(Z_{s},~s=(t,i)\in\mathcal{S})$, a spatial field on
$\mathcal{S}=\mathcal{T}\times\mathcal{I}$ where $\mathcal{T}$ is defined as
previously and $\mathcal{I}=\{1,2,\cdots,m\}$; we consider a finite space state $F$, and for each $s\in\mathcal{S}%
,~Z_{s}\in F$; without loss of generality, we assume
that the distribution $\pi$ of $Z$ is a Gibbs distribution with translation
invariant potentials $\Phi_{A_{k}}(\bullet),~k=1,K$ associated with a family
$\mathcal{A}=\{A_{k},k=1,K\}$ of subsets of $\mathcal{S}$, and that $\Phi_{A_{k}}(z)$
depend only on $z_{A_{k}}$, the layout of $z$ over $A_{k}$. Then $\pi$ is
characterized by the energy:%
\[
U(z)=\sum_{k=1,K}\ \sum_{s\in S(k)}\Phi_{A_{k}+s}(z),
\]
with $S(k)=\{s\in\mathcal{S}\ \mathrm{\ s.t.}\ A_{k}+s\subseteq\mathcal{S}\}.$

The basic idea is to consider the spatial Gibbs field with states in $F$ as a
$m$-multivariate Gibbs process on $E=F^{m}$; with this intention, we are going
to write the energy as a sum of time potentials.

We write $Z_{t}=(Z_{(t,i)},~i\in\mathcal{I})\ (Z_{t}\in E=F^{m})$ and adopt
the notation $z_{t}=(z_{(t,i)},~i\in\mathcal{I})$.

For $A\subseteq\mathcal{S}$, we define the height of $A$ by $r(A)=\sup
\{\left\vert u-v\right\vert ,\exists(u,i)$ and $(v,j)\in A\},$ and we define
$r=r(\mathcal{A})=\sup\{r(A_{k}),~k=1,K\}$ to be the largest height of the
potentials. With this notation, we write again the energy $U$
\[
U(z)=\sum_{h=0}^{r}\sum_{t=h+1}^{T}\Psi(z_{t-h},\cdots,z_{t})
\]
with $\Psi(z_{t-h},\cdots,z_{t})=\sum_{k:r(A_{k})=h}\ \sum_{s\in S_{t}(k)}%
\Phi_{A_{k}+s}(z)$ where $S_{t}(k)=\{s=(u,i):A_{k}+s\subseteq\mathcal{S}$ and
$t-r(A_{k})\leq u\leq t\}$.

Then $(Z_{t})$ is a Markov random field w.r.t. the $2r$-nearest neighbors system
but also a Markov process with memory $r$: $Y_{t}=(Z_{t},Z_{t+1},\cdots
Z_{t+r}),$ $t=1,T-r$, is a Markov chain on $E^{\ast}=E^{r}$ for which we get
the results (\ref{lois marginales}) and (\ref{calcul de la constante}).

We give a concrete example of the procedure below.

\subsection{Computing the normalization constant for an Ising model}

\label{Ising spatial} Let us illustrate the result with a translation invariant Ising
model. Let $\mathcal{S}=\mathcal{T}\times\mathcal{I=}\{1,2,\cdots
,T\}\times\{1,2,\cdots,m\}$ be the set of sites, $F=\{-1,+1\}$ the state
space, and $Z=(Z_{(t,i)},~(t,i)\in\mathcal{S)}$ the Markov random field on
$\mathcal{S}$ with the $4$-nearest neighbors system. The joint distribution
$\pi$ of $Z$ is characterized by the potentials
\begin{align*}
\Phi_{t,i}(z)  &  =\alpha~z_{(t,i)}\ \ \ \mathrm{\ for}\ (t,i)\in
\mathcal{S},\\
\Phi_{\{(t,i),(t,i+1)\}}(z)  &  =\beta~z_{(t,i)}z_{(t,i+1)}\ \mathrm{\ for}%
\ 1\leq i\leq m-1,\\
\ \mathrm{and}\ \Phi_{\{(t,i),(t+1,i)\}}(z)  &  =\delta~z_{(t,i)}%
z_{(t+1,i)}\ \mathrm{\ for}\ 2\leq t\leq T.
\end{align*}
We write $Z$ as a temporal Gibbs process with the following potentials:
\begin{align*}
\Psi(z_{t})=\alpha\sum_{i=1,m}z_{(t,i)}+\beta\sum_{i=1,m-1}z_{(t,i)}%
z_{(t,i+1)},\\
\Psi(z_{t},z_{t+1})=\delta\sum_{i=1,m}z_{(t,i)}z_{(t+1,i)},\ 1\leq t\leq T-1.
\end{align*}
Denoting the cardinal of $G$ as $\sharp G$, we define the following counting statistics associated with $c$, $d\in
E=\{-1,+1\}^{m}$: $n^{+}(c)=\sharp\{i\in\mathcal{I}:c_{i}=+1\}$,
$n^{-}(c)=m-n^{+}(c)$, $v^{+}(c)=\sharp\{i=1,m-1:c_{i}=c_{i+1}\}$,
$v^{-}(c)=(m-1)-v^{+}(c)$, and finally $n^{+}(c,d)=\sharp\{i\in\mathcal{I}%
:c_{i}=d_{i}\}$, $n^{-}(c,d)=m-n^{+}(c,d)$.

We then apply formula (\ref{calcul de la constante}) w.r.t. the following
$2^{m}\times2^{m}$ matrices $(H_{t})$ defined for any $u,v\in E=\{-1,+1\}^{m}$
and $\delta(t)=\delta\times\mathbf{1}_{(t\leq T-1)}$ to be:

$H_{t}(u,v)=\exp\left\{\alpha(n^{+}(u)-n^{-}(u))+\beta(v^{+}(u)-v^{-}%
(u))+\delta(t)(n^{+}(u,v)-n^{-}(u,v))\right\}.$

\bigskip

The limitation of the method is linked to the dimension $m$ of
the lattice since we deal with $2^{m}\times2^{m}$ matrices; on the other hand,
there is no constraint on $T$, the ``time'' dimension of the lattice.

\bigskip
\textit{Example 4: Ising model on a lattice}

We consider the Ising model described above on a lattice $m\times T$ with
$m=10$, and we compute $C$ for increasing values of $T$ using
(\ref{calcul de la constante}) and for the two methods, calculating the power
$H^{T-2}$ of $H $, or using a diagonalization of $H$.

\begin{table}
\caption{Computing times of $C$ for an Ising model on a lattice $10\times T$.}
\label{tab:4}       
\begin{tabular}
[c]{|c|c|c|c|}\hline
$m=10$ & M1 & M2 & $C$\\\hline
$T=10$ & \multicolumn{1}{|c|}{0.7907} & \multicolumn{1}{|c|}{3.4044} &
5.4083e+030\\\hline
$T=100$ & \multicolumn{1}{|c|}{1.2494} & \multicolumn{1}{|c|}{3.8640} &
2.4344e+307\\\hline
$T=500$ & 3.2169 & 5.7370 & 1.3868$\times2^{1023\times4}$e+305\\\hline
$T=1000$ & 3.3300 & 5.8851 & 2.1706$\times2^{1023\times9}$e+302\\\hline
$T=10000$ & 4.0314 & 6.5765 & 6.8918$\times2^{1023\times99}$e+251\\\hline
$T=10^{6}$ & 7.2580 & 9.7344 & 1.9362$\times2^{1020\times9999+600\times19}
$e+295\\\hline
\end{tabular}
\end{table}

We observe in Table \ref{tab:4} that it is computationally more efficient to compute the
powers $H^{T-2}$ rather than to use a diagonalization. Indeed, the Matlab
diagonalization procedure itself is expensive for large size matrices, and we
treat $2^{10}\times2^{10}$ matrices here. Moreover, the difference between the
execution times of the two methods is stable, around 2.5 seconds.

One might use other diagonalization algorithms, faster and computationally less expensive than the Matlab one. We give below a general idea of such a procedure.
Let us consider a $m\times n$ matrix $H$ where $m$ and $n$ are large numbers. We want to factor $H=PDQ^T$, with orthogonal matrices $P$ and $Q$ and a diagonal matrix $D$;  the algorithm is the following:

Step 1: Find a $m\times k$ orthogonal matrix $O$ such that $OO^TH$ is a good approximation of $H$. This step is performed with a probabilistic algorithm (see Halko 2011).

Step 2: Form the $k\times n$ matrix $B=O^TH$.

Step 3: Compute the decomposition of the small matrix $B=RDQ^T$.

Step 4: Form $P=OR$ which is still orthogonal.

We then obtain $H\simeq  PDQ^T$.

\subsection{Some generalizations}

The results can be extended to larger potentials, such as triplet potentials and so on.

We can also consider varying state spaces $E_{t}\ni z_{t}$;  recurrences
(\ref{recurrence fondamentale}) and (\ref{lois marginales}),
(\ref{calcul de la constante}) still hold with rectangular matrices $H_{t}$.

Finally, we can extend results for the marginals (\ref{recursion arriere}), the
normalization constant (\ref{calcul de la constante}) and recursive properties
(\ref{recurrence fondamentale}) to a sequence of embedded subsets of
$\mathcal{T}=\{1,2,\cdots,T\}$. Let us consider a decreasing sequence
$\mathcal{T}=S_{Q}\supset S_{Q-1}\supset\cdots\supset S_{1}$ of subsets of
$\mathcal{T}$ and assume $S_{q}=S_{q-1}\cup\partial S_{q-1}$ for $q=1,Q-1$.
Similarly to the previous future-conditional contributions
(\ref{lois marginales}) used in proposition \ref{new recursion}, we define the
contributions as $\gamma_{q}(z(S_{q});z(\partial S_{q}))$, the future-conditional
energy as $U_{q}^{\ast}(z(S_{q});z(\partial S_{q}))=U_{q-1}^{\ast}(z(S_{q-1}%
);z(\partial S_{q-1}))+\Delta_{q}(z(\partial S_{q-1}));z(\partial S_{q})),$
with $\Delta_{q}(z(\partial S_{q-1}));z(\partial S_{q}))=\displaystyle\sum
_{u\in\partial S_{q-1}}\theta_{u}(z_{u})\ +\sum_{u\in\partial S_{q-1}%
,v\in\partial S_{q},<u,v>}\Psi_{\{u,v\}}(z_{u},z_{v}),$ and then the matrices
$H_{q}$ as: $ H_{q}(\partial S_{q};\partial S_{q-1})=\exp\Delta_{q}(z(\partial
S_{q-1});z(\partial S_{q}))$.

\bigskip

As an illustration let us set the following decreasing sequence $\mathcal{T}%
=S_{T-1}=\{1,2,\cdots,T\}$, $S_{T-2}=\{1,2,\cdots,T-1\}$, $\cdots$,
$S_{2}=\{1,2,3\}$ and $S_{1}=\{2\}.\ $\ \ For $q=T-1,\cdots,3,$ the
conditional contributions and the matrices $H$ are defined in the usual way
with (\ref{future cond energy}) and $H_{q}(u,v)=\exp\{\theta_{q}(u)+\Psi
_{q}(u,v)\}$, while for $q=2,~H_{2}$ is a $N\times N^{2}$ matrix with
$H_{2}(z_{4},(z_{1},z_{3}))=\exp\left\{\theta_{1}(z_{1})+\theta_{3}(z_{3})+\Psi
_{4}(z_{3},z_{4})\right\}.$

\section{Conclusion}

This paper gives results permitting the evaluation of any marginals and the
normalizing constant for temporal Gibbs processes $\pi$. We also give some
example for which, by a direct use of the Markovian properties, we can carry
maximum likelihood estimation without the normalizing constant. We provide the
same results for Gibbs random fields, treating the field as a temporal
multidimensional process.

Computing the normalizing constant as a matrix product is easy to implement
and fast; to increase the speed of the computation, it may useful to
diagonalize first; we suggest comparing for some high power of a matrix $H$ with
 diagonalization computing times in order to choose the method. Moreover,
in the spatial $T\times m$ case, while we can increase the temporal parameter $T$ without constraints,  the limitation of the procedure arises from $m$ due to the
manipulation of $N^{m}\times N^{m}$ matrices. So the method seems to fail for
large square lattices. As a comparison, Pettitt et al. (2003) compute the
normalizing constant for an auto-logistic model defined on a cylinder lattice
for which the smallest row or column is not greater than 10. They suggest to
split large lattices into smaller sublattices along the smallest row or
column. Similarly, Friel and Rue (2007) deal with Markov fields on lattices of
size up to $m\leq19$ (depending on other parameters); for larger lattices,
they partition a $46\times72$ lattice into 3 disjoint sublattices of dimension
$46\times17$ and a final sublattice of dimension $46\times18$, each separated
by a column of length 46. Similar ideas could apply here.
Another solution is to use directly the matrix formulation to overcome the curse of dimensionality. For instance, low-rank matrix approximation replaces the large-scale $m\times n$ matrix $A$ by a product $BC$ of respectively $m\times k$ and $k\times n$ matrices $B$ and $C$ with $k<m\wedge n$ (see Witten 2013); $B$ is a subset of columns of $A$.  This will speed-up  matrix-vector and matrix-matrix multiplication. We can also use randomized (approximate) methods to find an eigenvalue decomposition, as described in the previous section. These examples show that the matrix formulation opens the possibility of fast computation methods.

\end{document}